\documentclass[12pt]{article}
\usepackage{epsfig,amsfonts, amssymb,amsmath, color, graphics,float}
\usepackage[colorlinks=true]{hyperref}

\newtheorem{prop}{Proposition}
\newtheorem{thrm}{Theorem}

\def \R{\mathbb R}
\def \no{\noindent }
\def \di{\displaystyle}

\def \var{\varepsilon}

\def \bp{{\bf Proof:}}
\def \ep{\hskip1.5cm {\tiny {$\blacksquare$}}\\}

\title{Study of a degenerate reaction-diffusion system arising in particle dynamics with aggregation effects}
\author{Mich\`ele GRILLOT$^1$, Philippe GRILLOT$^1$ and Simona MANCINI $^1$}
\date{ }

\begin{document}
\maketitle

\noindent
{\bf Keys words:}  degenerate parabolic system, existence of solutions, stability, pattern formation, aggregation\\


\addtocounter{footnote}{1} \footnotetext{ Math\'ematiques - Analyse, Probabilit\'es, Mod\'elisation - Orl\'eans (MAPMO),  CNRS, UMR 7349  Universit\'e d'Orl\'eans - Route de Chartres - B.P. 6759 - 45067  Orl\' eans cedex 2 - France \\ 
michele.grillot@univ-orleans.fr - philippe.grillot@univ-orleans.fr - simona.mancini@univ-orleans.fr}

\abstract{In this paper we are interested in a degenerate parabolic system of reaction-diffusion equations arising in biology when studying cell adhesion at the protein level. In this modeling the unknown is the couple of the distribution laws of the freely diffusing proteins and of the fixed ones. Under sufficient conditions on the aggregation and unbinding probabilities, we prove the existence of solution of the considered system, as well as their positivity, boundedness and uniqueness. Moreover, we discuss the stability of the equilibrium solution. Finally, we show that the simpler and particular choice of an affine aggregation function and of a constant unbinding probability do not lead to pattern formation as expected in the application. These analytical results are also supported by some numerical simulation.\\}

\section{Introduction }\label{intro}

The study and modeling of the adhesion process between cells is an active and complex research subject both in biology and in mathematics, and has implications, for example, with the cancer growth and the formation of neuronal connections in the embryon (see for example \cite{bio1}, \cite{bio2}, \cite{bio3},\cite{math-bio1}, \cite{math-bio2} and \cite{math-bio3}). The main protein involved in the growth of adhesion junctions between cells is called cadherin. Cadherins diffuse on the cell membrane and
may link to other cadherins on the same membrane (cis-contact) as well as bind to other cadherins on another cell membrane (trans-contact). These localized contacts can then evolve and give rise to larger adhesion domain, called adhesion junctions. These adhesion junctions are linked to actin filaments and help to their polymerization and growth, which in turns leads to the cell migration. In the experimental framework we consider, the link with the actin filaments is inhibited, and therefore the adhesion junctions and the binding cadherins are fixed in space. More precisely, experimentally, a cell is left to expand on a microscope slide covered by a uniform substrate of fixed cadherin targets which mimic the membrane of a second cell. Thus, when a diffusing cadherin link to a target it stops diffusing and fixes on the position of the given target. It appears that, even without the influence of the actin filaments, contacts between cadherins are easier created next to already existent ones, this local effect creates aggregates of linked cadherins enforcing the cell-cell adhesion junctions. On the other hand, if the strength of the contact is not large enough, the initial link may break and the cadherin diffuses again. \\

From now on, we shall consider cadherins as particles and we describe their spatio-temporal dynamics on the substrate by a mathematical model at the mesoscopic scale. The unknown is the couple of density functions $(u,v)$ representing the particle distributions of the diffusing particles, $u=u(x,t)$, and of the fixed ones, $v=v(x,t)$, where $x\in \Omega \subset \R^2$ is the bi-dimensional space position and $t\geq 0$ is time, $\Omega$ is the bounded domain where the evolution takes place. We describe the evolution of the population of the diffusing particles, $u$, by a reaction-diffusion equation, and the evolution of the population of the fixed particles, $v$, by an non-linear ordinary differential equation, since the diffusion coefficient for $v$ is null. Hence we are dealing with a degenerate system of reaction-diffusion equations. The reaction term represents the change of state of a particle, form free to fixed and vice-versa. It is thus composed of two terms, a gain one describing the growth of the fixed particle population $v$, and a loss one giving the amount of broken links. The gain term is proportional to the number of free targets (i.e. those targets which have no particle linked to), to the density of diffusing particles and to the aggregation probability function describing the probability of a diffusing particle to link to a target. The loss term is proportional to the density of fixed particles and to a unbinding probability function. The mathematical analysis of this problem is done in a general framework, but in the sequel we prove that the simple definition of the aggregation function as an affine function of $v$, and of the unbinding function as constant, is not in agreement with the experimental observations, since there is no pattern formations, so that no isolated adhesion junctions can be obtained. \\

In this paper, we are mainly interested by the theoretical study of the solutions $(u,v)$. The main results we prove are the existence and uniqueness of solutions to the degenerate reaction-diffusion system and their positivity and boundedness. We also discuss the existence of solutions to the stationary problem and their stability in other to prove the existence of pattern formations.  All these results are interesting in order to show that the proposed model makes sense in the biological framework and that it is, or not, consistent with the experimental results. For instance, we show that the simple affine aggregation function and the constant unbinding probability are not appropriate to recover the results shown by experiments. In fact, for this simple choice, the stationary solution turns out to be homogeneous in space and stable with respect to small perturbations, so that no pattern formation is possible, which would correspond in our model to no formation of isolated particle aggregates. A non-linear reaction term based on a more detailed description of the aggregation and unbinding phenomena, together with the comparison with biological experimental results, is given in \cite{SMSM}. Although a large number of parabolic, degenerate systems or not, were already studied (see for exemple \cite{amann1, amann2, mar, rw, soup} and as well as the references therein), their results do not directly apply to our model because of the particular form of our reaction term. Note, in particular, that our reaction term behaves like a polynomial which doesn't satisfy the bound of \cite{amann2}.\\

The paper is organized as follows. In section \ref{model} we detail the mathematical  model describing the particles dynamics and we announce the main results. Section \ref{proof} is devoted to the proofs. In section \ref{stability} we study the stability of the stationary solution to the problem. Finally, in section \ref{application}, we 
consider the particular case of an affine aggregation function and constant unbinding probability, and show that the stationary solution in this case is stable. This result is also  illustrated by some numerical simulation.

\section{ The model and main result}\label{model}

We consider a bounded regular bi-dimensional space $\Omega \subset \mathbb R^2$.  Then, the unknown of our problem is the couple of density functions $(u,v)=(u(t,x), v(t,x))$ of both the diffusing and fixed particles. For instance, $u=u(t,x)$ (respectively $v=v(t,x)$) represents the probability to find a diffusing (respectively fixed) particle at time $t\geq 0$ in the position $x\in \Omega$. Note that $x$ is a two dimensional vector point. We assume the targets density on the substrate to be constant in time and space, and we denote it by $\rho$ so that, up to a normalisation, $0 < \rho \leq  1$. From a biological point of view, it make sens to assume that there can't be more fixed particles than targets.  Moreover, since $u$ and $v$ are density distributions, they must be non-negative. Thus it must be, for all time $t\geq 0$ and all $x \in \Omega$ :
\begin{equation}\label{un}   
0 \leq u(x) \quad \hbox{ and } \quad 0\leq v(x)\leq \rho .
\end{equation}
Note that $v(t,x)=\rho$ holds if and only if we have the critical situation in which fixed particles never unlink from the target.\\
Recalling that, free particles diffuse and become fixed when linking to a target, and that 
fixed ones don't diffuse and become free when unbinding from a target, the respective density  functions $u$ and $v$ must then satisfy the following degenerate system of reaction-diffusion equations :
\begin{equation}\label{deux} 
\begin{cases} 
\partial_t u-\sigma \Delta u = - Q(u,v)\\
\partial_tv=Q(u,v)
\end{cases} 
\end{equation}
in $(0,+\infty)\times \Omega$, where $\sigma >0$ is the diffusion coefficient. \\

In \eqref{deux}, the function $Q(u,v)$ is the reaction term which describes how particles pass from a state to the other one: diffusing to fixed and vice-versa. As explained previously, see section \ref{intro}, the reaction term $Q$ is composed by the difference of a gain term $Q^+$,  counting the diffusing particles that became fixed, and of a loss term $Q^-$, counting the fixed particles that became free. The gain term $Q^+$ is thus proportional to the density of free targets $\rho-v$, to the density of free particles $u$ and to an  aggregation probability function $F=F(v)$.
Analogously,  the loss term $Q^-$ is assumed to be proportional to the density of fixed particles $v$, and to the unbinding probability $G=G(v)$. 
Hence $Q$ is given by:
\begin{equation}\label{trois} 
Q(u,v)=Q^+(u,v) - Q^- (u,v) = u (\rho-v) F(v) - v G(v)
\end{equation}
Since \eqref{un} should hold, let us consider $v\in [0,\rho]$, and fix the assumption we need on $F$ and $G$ as follows: 
\begin{align*}
& a)\  \forall\ v\in [0,\rho],\quad 0< F(v), G(v) \leq 1 \\
& b)\ F\textrm{\ and\ } G \textrm{ \ are\  locally\ Lipschitz}
\end{align*}
Assumption a) reflects the fact that $F$ and $G$ are probabilities. Assumption b) is the minimal regularity assumption needed in the sequel. 
When needed we'll continuously prolong the definition of $F$ and $G$ to the whole real values $v\in \mathbb R$ just defining:
\begin{align*}
& c)\ \forall\ v \leq 0, \quad F(v)=F(0) \; , \quad G(v)=G(0)\\
& d)\ \forall\ v\geq \rho, \quad F(v)=F(\rho)\; , \quad G(v)=G(\rho)
\end{align*}

Note that, in the biological framework, the aggregation function $F$ must be increasing in $v$, but  this condition is not required for the existence positivity and boundedness proofs. In Section \ref{application} we'll consider a particular case in which $F$ is a strictly increasing affine function for $v\in [0,\rho]$ and $G$ is constant (see \eqref{Fdef} and \eqref{Gdef}). With these choices $F$ and $G$ satisfy assumptions a) and b), but also relation \eqref{relation-h(v)}, leading to a unique homogeneous stable solution at equilibrium.\\

Since $\Omega$ is bounded, the first equation in \eqref{deux} needs to be endowed by some boundary conditions. We choose to describe them by Neumann boundary conditions:
\begin{equation}\label{cinq}  
{{\partial u}\over{\partial \nu}}(t,x)=0 ,
\end{equation}
biologically implying a zero flux on the boundary $\partial \Omega$.

Further, we complete system \eqref{deux} by the initial condition, for $x \in \Omega$, $(u_0,v_0)=(u(0,x), v(0,x))$ representing the initial distribution of both diffusing and fixed particles. 
We impose $u_0$ and $v_0$ continuous in $\overline{\Omega}$ and such that \eqref{un} holds $\forall\ x\in \Omega $, that is:
\begin{equation}\label{hypo1}
0 \leq u_0(x) \quad \textrm{and} \quad 0 \leq v_0(x) \leq \rho.
\end{equation}
Since $u_0$ is a continuous function on a bounded domain $\Omega$, we define its maximum on $\Omega$ as: 
\begin{equation}\label{max} \mu = \max_{x\in \Omega} u_0(x). \end{equation} 

Up to a multiplicative factor, we can assume the total initial density of particles to be  normalized to 1, that is :
\begin{equation}\label{sixbis}
\int_{\Omega} u_0(x)+v_0(x) \ dx =1
\end{equation}
Since particles aren't loss or created, but just change their status from free to fixed and vice-versa, \eqref{sixbis} must hold at any time $t\geq 0$, that is the following conservation property must hold:
\begin{equation}\label{six}
\int_{\Omega} u(t,x)+v(t,x) \ dx =1 , \qquad \forall \ t \geq 0 .
\end{equation}
\begin{prop} The solution of system \eqref{deux} completed by the boundary conditions \eqref{cinq} and with an initial data satisfying \eqref{sixbis}, satisfy the conservation property \eqref{six}.
\end{prop}
\bp \
Summing the equations in \eqref{deux}, integrating over $\Omega$, and applying the Green formula, we get :
$$ {{\partial\over{\partial t}} }\int_{\Omega} (u+v) dx -\sigma \int_{\partial\Omega} {{\partial u}\over{\partial \nu}}d\sigma =0,$$
and thanks to the boundary conditions \eqref{cinq} we obtain:
$$\frac{\partial}{\partial t} \int_\Omega (u+v) dx =0,$$
concluding the proof.
\ep

Following \cite{BS}, by a classical solution $(u,v)$ of \eqref{deux}, we mean $(u,v) \in (C([0,T)\times \overline{\Omega}))^2 \cap (C^{1,2}((0,T)\times  \overline{\Omega}) \times C^{1,0}((0,T)\times  \overline{\Omega}))$ solution of \eqref{deux}, where $T \in (0,\infty ]$ denotes its maximal existence time. Before enter in the details of the analysis of the solutions to system \eqref{deux}, we state here the main results, which will be proved in the next sections. We start by a result important from a biological point of view, since it ensures that the solutions to \eqref{deux} remains non-negative and bounded. 
\begin{thrm}\label{positivity}
Assume that the initial data $(u_0,v_0)$ satisfies \eqref{hypo1}, 
then a classical solution $(u,v)$ to \eqref{deux} verifies $\forall\ t>0$, $ \forall\ x\in \Omega$:
\begin{equation}\label{borne}
0 \leq u(t,x) \leq \mu +\rho t\ , \qquad 0 \leq v(t,x) \leq \rho .
\end{equation}
\end{thrm}
We express now our main result , which is also important since it ensures the existence of solutions to system \eqref{deux}.

\begin{thrm}\label{main-theorem}
Let $T>0$, let $F$ and $G$ satisfy conditions a) to d) and consider an initial data with both $u_0=u_0(x)$ and $v_0=v_0(x)$ in $C^0(\overline{\Omega})$ and satisfying \eqref{hypo1} for all $x\in \overline{\Omega}$. Then there exists two non-negative functions $ u\in W^{1,p}( (0,T)\times \Omega)$, $p>1$ and $v\in L^{\infty}((0,T)\times \Omega)$, solutions to \eqref{deux}-\eqref{cinq} endowed by the initial conditions: 
\begin{equation}\label{onze} 
 u(0,.)=u_0 \quad \hbox{ and } \qquad v(0,.)=v_0  
\end{equation}
\end{thrm}
Note that we don't have the unicity in Theorem \ref{main-theorem}.
If we add some monotonicity conditions on $F$ and $G$, then we have the following:
\begin{thrm}\label{uniq-theorem}
Let $T>0$, let $F$ and $G$ satisfy conditions a) to d) and be such that $v \mapsto -u(\rho - v)F(v)+vG(v)$ is nondecreasing. Let us consider an initial data $(u_0,v_0)$ such that both $u_0$ and $v_0$ are in $C^0(\overline{\Omega})$ and satisfying \eqref{hypo1}
for all $x\in \overline{\Omega}$. Then there exists a unique solution $(u,v)$ of \eqref{deux}-\eqref{cinq} and \eqref{onze}, with $v(t,x)=0 \ \forall (t,x) \in (0,T)\times \partial \Omega$ and such that $u \in C([0,T)\times \overline{\Omega})\cap C^{1,2}((0,T)\times  \overline{\Omega})$ and $v \in  C([0,T)\times \overline{\Omega})\cap C^{1,0}((0,T)\times  \overline{\Omega})$ are non-negative functions.
\end{thrm}

Moreover, note that both Theorem \ref{positivity} and Theorem \ref{uniq-theorem} imply that the solutions exist on $(0, \infty )\times \Omega$.\\

Once proved that the problem is well posed, we will consider (see section \ref{stability}) the existence and stability of solutions at equilibrium, that is of solution of the stationary problem associated to \eqref{deux}-\eqref{cinq}. 
We will prove that there exist at least one stationary solution which is homogeneous in space and satisfying \eqref{un}. Moreover, under the condition:
\begin{equation}\label{relation-h(v)}
\rho F(v)G(v) + v (\rho-v)(F(v) G'(v) - F'(v) G(v)) >0 ,
\end{equation}
we obtain uniqueness of the stationary homogeneous solution, which we denote by $(U,V)$, and also its stability. \\

Finally, defining $F(v)$ by \eqref{Fdef}, and $G(v)$ by \eqref{Gdef}, we will show that , the stationary solution $(U,V)$ is stable with respect to small perturbations, and thus we can't have formation of aggregate of particles. This can be overcome choosing a non-linear aggregation function and also a non-linear unbinding rate, as it is done in \cite{SMSM}.

\section{Existence, uniqueness and boundedness}\label{proof}

We consider first the positivity and boundedness of the solution $(u,v)$ to \eqref{deux}, assuming that it exists. We want thus to prove Theorem \ref{positivity}.

\subsection{Proof of Theorem \ref{positivity}}
We shall detail the proof in three steps. First we prove the upper bound for $v$.
Second we prove the bound for $u$. Third we prove the positivity of both $u$ and $v$.\\

\noindent
\bp\\
\underline{Step I.} 
We first prove by contradiction the upper bound for $v$, $v\leq \rho$. Let us define the set
$$A=\{t > 0\ :\ \exists\ x \in \Omega\ \textrm{such that}\ v(t,x) > \rho \},$$
and assume that there exists a time $t_0' >0$ and a point $x_0' \in \Omega$ such that $v(t_0,x_0') > \rho$. 
Hence, $t_0' \in A$ and $A$ is bounded from below by 0, so that it exists $t_0=\inf A$. \\
Note that, by means of continuity arguments, for this time $t_0$, it exists a point $x_0 \in \Omega$ such that $v(t_0,x_0)= \rho$. \\
Now, considering the equation for $v$ in \eqref{deux} and computing its value at the point $(t_0,x_0)$, we get, because of assumption a),
$$\partial_t v(t_0,x_0) = -G(\rho) \rho <0.$$
Therefore, the function $v(t,x_0)$ is decreasing in time in a neighborhood of $t_0$.
Hence, it exists $\eta >0$ such that for all $t\in [t_0-\eta, t_0]$ we have $v(t,x_0)>v(t_0,x_0)=\rho$. Thus $t_0-\eta \in A$ and $t_0-\eta <t_0$, contradicting 
the fact that $t_0=\inf A$. Hence, $t_0'$ doesn't exists and $A=\emptyset$.\\

\noindent
\underline{Step II.} 
Having that for all time $t\geq 0$ and all $x \in \Omega$, $v\leq \rho$, we now consider the upper bound for $u$, $u\leq \mu + \rho t$. \\
We prove that the function 
$(u-\mu-\rho t)^+ =\max(0, u-\mu - \rho t) $
is equal to zero for all time $t\geq 0$ and all $x\in \Omega$, where $\mu$ is defined in \eqref{max}. \\
Note that if $(u-\mu-\rho t)^+ \neq 0$, then $u\geq \mu+\rho t $, so that $u >0$. Using the equation for $u$ in \eqref{deux}, and since, from a), $G(v) \leq 1$ and, from Step I, $v\leq \rho$, we have :
\begin{align*}
\partial_t (u-\mu-\rho t) &(u-\mu-\rho t)^+ -\sigma \Delta (u-\mu-\rho t)(u-\mu-\rho t)^+\\
&=(\partial_t u -\sigma \Delta u)(u-\mu-\rho t)^+ - \rho (u-\mu-\rho t)^+\\
&=-u(\rho-v)F(v)(u-\mu-\rho t)^+ + (vG(v) -\rho)(u-\mu-\rho t)^+ \leq 0 .
\end{align*}
Hence, integrating over $\Omega$ :
$$\frac{1}{2}\frac{d}{dt} \int_\Omega ((u-\mu-\rho t)^+)^2 \ dx + \sigma \int_\Omega 
| \nabla (u-\mu-\rho t)^+|^2\ dx \leq 0,$$
so that, for all $t\geq 0$, it must be 
$$\| (u-\mu-\rho t)^+ \|_2^2 \leq  \| (u(0,x)-\mu)^+\|_2^2 = 0.$$
Concluding this second step of the proof.\\

\noindent
\underline{Step III.} Let $(u,v)$ be a solution of \eqref{deux}, under hypothesis \eqref{hypo1} we prove that $v \geq 0$ and $u \geq 0$ for all times $t \leq T$. Let $T' \in (0,T)$. Multiplying equations of \eqref{deux} respectively by $u^-=\min (u,0)$ and $v^-=\min(v,0)$, integrating on $\Omega$ and using Green's formula, we obtain:
\begin{align*}
\displaystyle {1 \over 2} {d \over {dt} }&\int_{\Omega}[(u^-)^2+(v^-)^2] dx + \sigma \int_{\Omega} | \nabla(u^-)|^2 dx\\
\displaystyle  &=-\int_\Omega (u^-)^2(\rho -v)F(v) dx +\int_\Omega v^+G(v)u^- dx +\int_\Omega v^-G(v)u^- dx  \\
\displaystyle &+\int_\Omega u^+v^-(\rho -v)F(v) dx+\int_\Omega u^-v^-(\rho -v)F(v) dx -\int_\Omega (v^-)^2G(v) dx
\end{align*}
where we have used that $u=u^++u^-$ and $v=v^++v^-$. Because of a), c) and d), because of Step I, and since $u^+v^- \leq 0$ and $u^-v^+ \leq 0$, we get:
\begin{align*}
\displaystyle {1 \over 2} {d \over {dt} }&\int_{\Omega}[(u^-)^2+(v^-)^2] dx \\
&\leq \int_\Omega (\rho -v)F(v) [- (u^-)^2 +u^-v^-] dx + \int_\Omega G(v) [-(v^-)^2 +u^-v^-]  dx.
\end{align*}
But $u^-v^-\leq [(u^-)^2+((v^-)^2]$, implies:
\begin{equation}\label{mille} 
 {d \over {dt} }\int_{\Omega}[(u^-)^2+(v^-)^2] dx \leq \int_\Omega (\rho -v)F(v) (v^-)^2 dx + \int_\Omega G(v) (u^-)^2 dx.
\end{equation}
Recalling that $v \in C([0,T)\times \overline{\Omega})$, we denote by $m=m(v,T')$ the minimum of $v$ 
in $[0,T']\times \overline{\Omega}$, and using the upper bound in a), we deduce from \eqref{mille} that:
$${d \over {dt} }\int_{\Omega}[(u^-)^2+(v^-)^2] dx \leq \max (\rho -m, 1)  \int_{\Omega}[(u^-)^2+(v^-)^2] dx.$$
Since $u_0 \geq 0$ and $v_0 \geq 0$, the Gronwall's lemma allows us to conclude that
$$\int_{\Omega}[(u^-)^2+(v^-)^2]\leq 0,$$
so that nor $u$ nor $v$ can be negative on $[0,T]\times \overline{\Omega}$ . \ep

\subsection{Proof of Theorem \ref{main-theorem} and Theorem \ref{uniq-theorem}}

We prove here the existence of solution to problem \eqref{deux}-\eqref{cinq} and \eqref{onze}, by constructing two Cauchy sequences $(u_n)_{n\geq 0}$ and $(v_n)_{n\geq 0}$ converging to $u$ and $v$ solution to \eqref{deux}-\eqref{cinq} and \eqref{onze}.\\

 Let $T>0$ and let $u_0$ and $v_0$ be two continuous function on $\overline{\Omega}$ and satisfying \eqref{hypo1}. For all $(t,x)\in [0,T]\times \overline{\Omega}$, let us define two sequences $(u_n)_{n\geq 0}$ and $(v_n)_{n\geq 0}$ by:
\begin{equation}\label{quatorze} 
\begin{cases} 
\displaystyle \partial_t u_{n+1} -\sigma \Delta u_{n+1} =-Q(u_{n+1},v_{n}) \quad \hbox{ in } (0,T)\times \Omega \\
\displaystyle \partial_t v_{n+1} =Q(u_n,v_{n+1}) \quad  \hbox{ in } (0,T)\times \Omega \\
\displaystyle {\di {{\partial u_{n+1}}\over{\partial \nu}} =0 \quad  \hbox{ on } (0,T)\times \partial \Omega}\\
\displaystyle u_{n+1}(0,.)=u_0  \quad \hbox{ in } \Omega \\
\displaystyle v_{n+1}(0,.)=v_0 \quad  \hbox{ in } \Omega ,
\end{cases}
\end{equation}
with \begin{equation}\label{treize}
u_0(t,x)=u_0(x)\quad \textrm{and} \qquad v_0(t,x)=v_0(x).
\end{equation}
Note that those functions well exist and are continuous, in fact the first equation is a linear heat equation in $u$ with a continuous coefficient $(\rho-v_n)F(v_n)$ 
and the second equation is a Riccati equation in time.\\
We'll prove that $(u_n)_{n\geq 0}$ and $(v_n)_{n\geq 0}$ are two Cauchy sequences, thus converging to two functions $u$ and $v$, which are solution to \eqref{deux}-\eqref{cinq} and\eqref{onze}. We first need some technical results.

\begin{prop}\label{propo1} 
Under the assumptions of the theorem both functions $u_n$ and $v_n$ defined by \eqref{treize} and \eqref{quatorze} satisfy the bound of Theorem \ref{positivity}, see \eqref{borne}, that is:
$$
0\leq u_n(t,x)\leq \mu + \rho t
$$
and
$$
0\leq v_n(t,x)\leq  \rho
$$
for all $(t,x)\in (0,T)\times \Omega$.
\end{prop}
\bp\\ 
The proof of this proposition is by induction using the same arguments as in the proof of Theorem \ref{positivity}. We omit it. \ep

The goal is to prove that both sequences ($u_n$) and ($v_n$) are Cauchy sequences. For simplicity of notations, for $n\in \mathbb N^*$ and $t\in (0,T)$, let us define by $U_n(t)$ and $V_n(t)$ the squares of the $L^2$-norms $\| u_{n+1}-u_n\|_{L^2(\Omega)}^2$ and  $\| v_{n+1}-v_n\|_{L^2(\Omega)}^2$, so that:
\begin{equation}\label{27}
U_n(t) =\int_{\Omega} (u_{n+1}-u_n)^2 dx \; ,\qquad V_n(t) =\int_{\Omega} (v_{n+1}-v_n)^2 dx.
\end{equation}
We first prove some technical results.
\begin{prop}\label{propo2} 
Under the assumptions of Theorem \ref{main-theorem}, there exists a constant $k>1$ such that for all $n\in \mathbb N^*$ and $t\in (0,T)$ :
\begin{equation}\label{28} 
U'_n(t) \leq 3k \, U_n(t) +k \, V_{n-1}(t) 
\end{equation}
and
\begin{equation}\label{29} 
V'_n(t) \leq 3k \, V_n(t) +k \, V_{n-1}(t). 
\end{equation}
\end{prop}
Note that \eqref{28} and \eqref{29} are completely analogous.\\

\noindent
{\bf Proof :}\\
The function $Q(r,s)$ defined by \eqref{trois} is a Lipschitz function in $[0,\gamma]\times[0,\rho]$, where 
$$\gamma=\mu + \rho T.$$ 
Therefore there exists a constant $k=k(\gamma,\rho)>1$ such that for all $((r,s),(r',s'))\in ([0,\gamma]\times[0,\rho])^2$ :
\begin{equation}\label{30} 
|Q(r,s)-Q(r',s')|\leq k(|r-r'|+|s-s'|).
\end{equation}
Since both functions $u_n$ and $u_{n+1}$ satisfy (\ref{quatorze}), we have:
\begin{equation}\label{31} 
\partial_t(u_{n+1}-u_n)-\sigma \Delta(u_{n+1}-u_n)=-Q(u_{n+1},v_n) +Q(u_n,v_{n-1})
\end{equation}
in $(0,T)\times \Omega$. Multiplying \eqref{31} by ($u_{n+1}-u_n$) and integrating on $\Omega$, we obtain for all $t\in (0,T)$:
\begin{equation}\label{32} 
{1\over 2} {d\over{dt}} \int_{\Omega} (u_{n+1}-u_n)^2 dx +\sigma \int_{\Omega} |\nabla(u_{n+1}-u_n)|^2 dx
\end{equation}
$$ \hskip2cm= \int_{\Omega} \left( Q(u_n,v_{n-1})-Q(u_{n+1},v_n)\right) \left( u_{n+1}-u_n\right) dx.$$
Because of \eqref{30} and \eqref{27} we deduce from \eqref{32} :
\begin{align*} 
{\di {1\over 2} \, U'_n(t)} &{\di  \leq k \int_{\Omega} \left( |u_{n+1}-u_n| +|v_n-v_{n-1}|\right) |u_{n+1}-u_n| dx }\\
& {\di  \leq k \int_{\Omega} |u_{n+1}-u_n|^2dx +{k\over 2} \int_{\Omega} \left( |v_n-v_{n-1}|^2 +|u_{n+1}-u_n|^2\right) dx} \\
&  {\di \leq {{3k}\over 2} U_n(t) +{k\over 2} V_{n-1}(t)}
\end{align*}
which implies \eqref{28}. The proof of \eqref{29} is completely similar because the second term of the first member of \eqref{32} doesn't play any role. \ep

The differential inequalities of Proposition \ref{propo2} allow us to obtain over-bounds inequalities of $U_n$ and $V_n$ :

\begin{prop}\label{propo3} 
Under the assumptions of Theorem \ref{main-theorem}, we have for all $n \in \mathbb N^*$ and $t\in (0,T)$ :
\begin{equation}\label{34} 
U_n(t) \leq kL\left( ke^{3kT}\right)^n {{t^n}\over{n!}}
\end{equation}
and similarly :
\begin{equation}\label{35} 
V_n(t) \leq kL \left( ke^{3kT}\right)^n {{t^n}\over{n!}}
\end{equation}
where $k\geq 1$ is the Lipschitz constant of  $Q$ in (\ref{30}) and $L$ is defined by :
$$
L=4[\max (\gamma,\rho)]^2 |\Omega |,
$$
where we denote by $|\Omega|$ the measure of $\Omega$ that is, $|\Omega|=\di \int_\Omega dx$.
\end{prop}
\bp\\
We prove it by induction. \\
For $n=1$, integrating \eqref{28} on $(0,T)$ for $t\in(0,T)$, we have:
$$
U_1(t)-U_1(0) \leq 3k \int_0^t U_1(s)dx +k\int_0^t V_0(s)ds .
$$
Recalling \eqref{quatorze} and \eqref{27}, it is easily seen that $U_1(0)= 0$. Note that by the same arguments we have $U_n(0)=0$ for all $n\in \mathbb N^*$. On the other hand, Proposition \ref{propo1} implies:
$$
V_0(s)=\int_{\Omega} (v_1(s,x)-v_0(s,x))^2dx \ \leq \int_{\Omega} (|v_1(s,x)|+|v_0(s,x)|)^2dx \ \leq 4\rho^2|\Omega|.
$$
Thus, for all $t\in(0,T)$ :
$$
U_1(t)\leq 4\rho^2 |\Omega| kt +3k\int_0^t U_1(s)ds.
$$
Gronwall's lemma leads us to the following inequality, for all $t\in (0,T)$ :
$$
U_1(t)\leq 4\rho^2|\Omega| kt +3k\int_0^t 4\rho^2 |\Omega| ks e^{3k(t-s)}ds .
$$
Integrating by parts we obtain for all $t\in(0,T)$ :
$$
U_1(t) \leq 4k\rho^2 |\Omega| \int_0^t e^{3k(t-s)}ds. 
$$
Since $k\geq 1$ and since $s\in [0,t]\subset [0,T]$, we have $t-s\leq T$ and:
$$
U_1(t)\leq 4k\rho^2 |\Omega| e^{3kT}t\leq Lk^2e^{3kT}t .
$$
In the same way we have for all $t\in(0,T)$ :
$$
V_1(t)\leq Lk^2e^{3kT}t.
$$
We prove now the induction step. Assume that \eqref{34} and \eqref{35} hold  for a fixed $n \geq1$, that is, for all $t\in(0,T)$:
$$\begin{cases}
\di U_n(t)\leq Lk^{n+1}e^{3kTn}{{t^n}\over{n!}}\\
\\ 
\di V_n(t)\leq Lk^{n+1}e^{3kTn}{{t^n}\over{n!}}
\end{cases}$$
From Proposition \ref{propo2}, since \eqref{28} and \eqref{35} hold, we have for all $t\in (0,T)$ :
$$
U'_{n+1}(t) \leq 3k \, U_{n+1}(t) +Lk^{n+2}e^{3kTn} {{t^n}\over{n!}}
$$
Therefore, integrating on $(0,t)$ and recalling that $U_{n+1}(0)=0$, we have:
$$
U_{n+1}(t)\leq 3k\int_0^t U_{n+1}(s) ds +Lk^{n+2}e^{3kTn} {{t^{n+1}}\over{(n+1)!}} .
$$
Gronwall's lemma implies for all $t\in (0,T)$ :
$$
U_{n+1}(t)\leq Lk^{n+2}e^{3kTn} {{t^{n+1}}\over{(n+1)!}} +3kLk^{n+2}e^{3kTn}\int_0^t {{s^{n+1}}\over{(n+1)!}} \, e^{3k(t-s)}ds.
$$
Integrating by parts we obtain for all $t\in (0,T)$:
$$
U_{n+1}(t) \leq Lk^{n+2}e^{3kTn}\int_0^t {{s^n}\over{n!}} \, e^{3k(t-s)}ds ,
$$
and since $t-s\leq t\leq T$, we deduce for all $t\in (0,T)$:
$$
U_{n+1}(t) \leq Lk^{n+2}e^{3kT(n+1)}{{t^{n+1}}\over{(n+1)!}}.
$$
In the same way we obtain 
$$V_{n+1}(t) \leq Lk^{n+2}e^{3kT(n+1)}{{t^{n+1}}\over{(n+1)!}} ,
$$
 which ends the proof. \ep

We now have all the ingredients to prove that $(u_n)_{n\geq 0}$ and $(v_n)_{n\geq 0}$  are Cauchy sequences.\\

\noindent
{\bf {Proof of Theorem \ref{main-theorem} :}\\} 
Let $t\in (0,T)$. We have for all $(n,m)\in (\mathbb N^*)^2$, $n>m$ :
$$
\| u_{n}(t,.)-u_m(t,.)\|_{L^2(\Omega)} \leq \sum_{j=m}^{n-1} \|u_{j+1}(t,.)-   u_{j} (t,.) \|_{L^2(\Omega)} \leq \sum_{j=m}^{n-1} (U_{j}(t))^{1/2}.
$$
Because of (\ref{34}) we deduce:
$$
\|u_{n}(t,.)-u_m(t,.)\|_{ L^2(\Omega)} \leq (kL)^{1/2} \sum_{j=m}^{n-1}  \left(\frac{(k e^{ 3kT} T)^j}{j !}\right)^{1/2}  .
$$
Since $\sum_{j=0}^\infty \left(\frac{C^j}{j !}\right)^{1/2}$, with $C>0$, is a convergent serie, we can conclude that $(u_n)_{n\geq 0}$ is a Cauchy sequence in $L^2((0,T)\times \Omega)$.
Similarly we prove that $(v_n)_{n\geq 0}$ is a Cauchy sequences in $L^2((0,T)\times \Omega)$. Therefore there exists two functions $u$ and $v$ in $L^2((0,T)\times \Omega)$ such that:
$$\|u_n-u\|_{L^2((0,T)\times\Omega)}\to 0\; , \quad\textrm{and}\qquad  
 \|v_n-v\|_{L^2((0,T)\times\Omega)}\to 0. $$
There also exists two subsequences, still denoted ($u_n$) and ($v_n$), respectively converging almost everywhere to $u$ and $v$ in $(0,T)\times \Omega$. Finally, Proposition \ref{propo1} implies:
$$0\leq u(t,x)\leq \mu+\rho t\; , \quad\textrm{and}\qquad
0\leq v(t,x)\leq \rho, $$
for almost all $(t,x)\in (0,T)\times \Omega$, so that both functions $u$ and $v$ are in $L^{\infty}((0,T)\times\Omega)$ and this is the same for $Q(u,v)$. Thus from classical  parabolic equations arguments, see \cite{lad}, we have that  $u\in W^{1,p}((0,T)\times\Omega)$,  $v \in L^{\infty}((0,T)\times \Omega)$ and they are solutions of \eqref{deux}-\eqref{cinq} and \eqref{onze}.\ep

Note that system \eqref{deux} doesn't satisfy classical comparison principle, and that the above theorems and proofs no more hold if $\rho$ is not a constant. In order to get uniqueness of the solution $(u,v)$ we need some more assumption on the aggregation function $F(v)$ and/or on the unbinding function $G(v)$. For instance, if $F(v)$ is nonincreasing and $G(v)$ nondecrasing, then the hypothesis of Theorem \ref{uniq-theorem} are satisfied. This will be the case for the particular choice \eqref{Fdef} and \eqref{Gdef}.\\

\noindent
{\bf Proof of Theorem \ref{uniq-theorem} :}\\ 
The conditions on $F$ and $G$ allow us to use the appendix of \cite{BS} and then we have a maximum principle for the problem which implies the uniqueness. The existence holds from theorem \ref{main-theorem} and from classical  parabolic equations arguments, see \cite{lad}. Note that because of the maximum principle the existence can also be proved using monotonous sequences.\ep

\section{Stability of solutions}\label{stability}
Now that the problem is proved to be well-posed, we consider its behavior at equilibrium, and deal, in this section, with the stability of the solution $(u,v)$ to \eqref{deux}-\eqref{cinq} completed by an initial data $(u_0,v_0)$ satisfying \eqref{onze}. The first step is to determine the solutions to the stationary problem associated to \eqref{deux}-\eqref{cinq}, which reads, for $x \in \Omega$:
\begin{equation}\label{sept}
\begin{cases} 
\Delta u=0\quad \textrm{in \ } \Omega\\
\displaystyle {\di {{\partial u}\over{\partial \nu}} =0 \quad  \hbox{ on } \partial \Omega}\\
Q(u,v)=0 \quad \textrm{in \ } \Omega.
\end{cases}
\end{equation}
These solutions are interesting since they give the solution at equilibrium of the considered problem. We have the following:
\begin{prop}\label{stationary-prop} 
Let $|\Omega |=1$, under assumptions a) to d), there exists at least one couple $(U,V)$ solution of problem \eqref{sept} satisfying the normalization condition \eqref{six} and bounds \eqref{un}. Such couple $(U,V)$ satisfies, for $V\in ]0,\rho[$: 
\begin{equation}\label{UV}
U= \frac{V G(V)}{(\rho-V)F(V)}
\end{equation}
Moreover, the unicity of the couple $(U,V)$ holds if the additional condition \eqref{relation-h(v)} is satisfied for all $v\in (0,\rho)$.
\end{prop}
{\bf Proof :}\\
\no Multiplying the first equation by $u$, integrating on $\Omega$ and applying the Green formula, we obtain that the solution $u$ of the first equation in \eqref{sept} must be homogeneous in $\Omega$, that is $u(x)=U\in \mathbb R$ for all $x\in \Omega$. Replacing $u=U$ in the last equation of \eqref{sept}, we then obtain that, if it exists, the solution $v$ of $Q(U,v)=0$ must be constant too, that is $v(x)=V\in \mathbb R$ for all $x\in \Omega$. From the relation $Q(U,V)=0$ we easily deduce \eqref{UV}.\\
We still have to prove the existence of such a solution satisfying \eqref{un} and \eqref{six}. From the normalization condition \eqref{six}, and since $|\Omega|=1$, we get that $U= 1 - V$. So that $V$ is implicitly defined by $Q(1-V,V)=0$ which leads to the relation :
\begin{equation}\label{relationV}
(1-V)(\rho-V)F(V)-VG(V)=0.
\end{equation}
We then seek for a value $V\in ]0,\rho[$ satisfying \eqref{relationV}.\\
 Let us define the following real function $h(v)$ for $v\in]0,\rho[$:
$$h(v)= \frac{vG(v)}{(\rho-v)F(v)}.$$
Then \eqref{relationV} reads:
$$h(v)+v-1=0.$$
Considering the limit values of the above functions at the boundaries $v=0$ and $v=\rho$, we have : $h(0)+0-1=-1 <0$ and $lim_{v\to \rho^-} h(v)+v-1 >0$. So by means of continuity arguments we can conclude that there exists at least one value $V\in ]0,\rho[$, such that \eqref{relationV} is satisfied. 
Moreover, since $U$ can be defined by \eqref{UV}, and thanks to assumption a) on $F$ and $G$, it easily follows that $U>0$. Thus both $U$ and $V$ satisfy \eqref{un}.\\
To have uniqueness for $V$, we need $h(v)$ to be strictly monotone. This is true if, for all $v\in ]0,\rho[$:
$$h'(v)=\frac{\rho F(v) G(v) + v (\rho -v) \left( F(v) G'(v) - F'(v) G(v)\right)}{(\rho-v)^2 F(v)^2}>0.$$ 
Thus, if \eqref{relation-h(v)} holds, $h(v)$ is strictly increasing in $(0,\rho)$, and so is the function $h(v)+v-1$, concluding the proof.\ep

It is well known that for some reaction-diffusion systems, some of the stationary solutions, or equilibrium points, for some ranges of the parameters, may not be stable, and that other non-homogeneous equilibrium structures may be found. This is known in morphogenesis as the Turing instabilities and happens in particular when the ratio of the diffusion coefficients is very small, see for exemple \cite{MurrayI, MurrayII}. In our problem, one of the diffusion coefficient is zero, and we are then interested to study the stability of the homogeneous stationary solution $(U,V)$, to see under which conditions patterns formation may occur. We have the following result:
\begin{prop}\label{stability-prop}
Under hypothesis a) to d), and if the homogeneous stationary solution $(U,V)$ satisfies \eqref{relation-h(v)}, then $(U,V)$ is a stable equilibrium solution of system \eqref{deux}-\eqref{cinq}.
\end{prop}
{\bf Proof :}\\
In order to study the stability of the equilibrium couple $(U,V)$, we must consider the linearized problem associated to \eqref{deux}-\eqref{cinq} around the equilibrium $(U,V)$. We need to evaluate the derivates of the reaction term $Q(u,v)$ with respect to $u$ and $v$ in $(U,V)$. Denoting these terms  respectively by :
$$A=\partial_u Q(U,V)\quad \textrm{and}\quad B=\partial_v Q(U,V),$$
and recalling that $U$ may be expressed as in \eqref{UV}, we get:
\begin{align*}
A&=(\rho-V)F(V),\\
B&=\frac{-V(\rho-V)
\left( F(V)G'(V)-F'(V)G(V)\right) -\rho F(V)G(V)}{(\rho-V)F(V)}.
\end{align*}
Since hypothesis a) to d) and \eqref{relation-h(v)} hold, then for $V\in ]0,\rho[$, $A>0$ and $B<0$ .\\
The linearized problem then reads :
\begin{equation}\label{linear}
\begin{cases}
\partial_t \tilde u =\sigma \ \Delta \tilde u - A \tilde u - B \tilde v\\
\partial_t \tilde v = A \tilde u + B \tilde v
\end{cases}
\end{equation}
where $\tilde u$ and $\tilde v$ are approximation of $u$ and $v$ such that $u=U+\tilde u$ and $v=V+\tilde v$. To solve this problem, under the Neumann boundary condition, we look to solution $(\tilde u, \tilde v)$ defined under the form :
\begin{equation}\label{fourier}
\tilde u = e^{\lambda t} \hat{u}(k), \quad \tilde v = e^{\lambda t} \hat{v}(k),
\end{equation}
where $\hat{u}(k)$ and $\hat{v}(k)$ are the Fourier transforms of $\tilde u$ and $\tilde v$ (i.e. $k>0$ is the Fourier wavenumber), and $\lambda$ are the eigenvalues which determines the temporal growth. It can easily be seen from \eqref{fourier} that, if $\lambda <0$, then both $\tilde u$ and $\tilde v$ goes to zero when $t$ goes to infinity, so that the solution $(u,v)$ converges back to the stationary homogeneous solution $(U,V)$, which is thus stable: a small perturbation of $(U,V)$ doesn't affect its stability.  We thus want to prove that the eigenvalues $\lambda$ are negative. Replacing \eqref{fourier} in \eqref{linear} and simplifying $e^{\lambda t}$, we get, for each $k$, the following linear system :
\begin{equation}\label{eigenvalues}
\begin{cases}
(\lambda +A+\sigma k^2) \hat u +  B \hat v = 0 \\
-A \hat u +(\lambda - B) \hat v=0.
\end{cases}
\end{equation}
Since we look for nontrivial solutions $(\hat u, \hat v)$ to \eqref{eigenvalues}, the determinant of the associated matrix must be zero, so that the values for $\lambda$ are determined has solutions of the characteristic polynomial associated to  \eqref{eigenvalues}. That is, for each $k$, $\lambda$ must solve :
$$\lambda^2 + \lambda(\sigma k^2 +A-B) - B\sigma k^2 =0.$$
That is,
\begin{equation}\label{lambda}
\lambda=\frac{1}{2} \left( -(\sigma k^2 +A-B) \pm \sqrt{(\sigma k^2 +A-B)^2 + 4B\sigma k^2} \right)
\end{equation}
Since $B<0$ and
$$(\sigma k^2+A-B)^2+4B\sigma k^2= A^2 + 2A\sigma k^2 - 2AB + (B+\sigma k^2)^2 >0,$$
there always exist two real solutions to \eqref{lambda}. Finally, since $(\sigma k^2+A-B) >0$,  in order to have $\lambda < 0$, we just have to check that the solution with the plus sign is negative. This holds if,
$$\sqrt{(\sigma k^2 +A-B)^2 + 4B\sigma k^2} < (\sigma k^2 +A-B),$$
which is true for all wavenumber $k>0$, since under hypothesis c), $B<0$, and 
conclude the proof. \ep \\

\section{A simple application}\label{application}
Getting back to the biology problem we inspired from, in this section we assume $F(v)$ to be an increasing and bounded affine function with respect to $v$, modeling the fact that the probability of a particle to be linked is larger when there already are fixed particles in its neighborhood :   
\begin{equation}\label{Fdef}
F(v) = \begin{cases}
\displaystyle a/(a+b\rho) \quad \textrm{if}\quad v<0\\
\displaystyle \frac{a + b v}{a+b\rho} \quad \textrm{if} \quad 0\leq v \leq \rho\\
1 \quad \textrm{if}\quad v>\rho
\end{cases}
\end{equation}
where $a >0$ represents the probability of a particle to naturally link to a target (i.e. if no other particle is fixed nearby) and $b >0$ is a parameter which control the influence of other fixed particles in the neighborhood. \\
Concerning $G(v)$, we choose it to be homogeneous with respect to $v$ : 
\begin{equation}\label{Gdef}
G(v)=\var\ , \quad \forall\ v \in \mathbb R,
\end{equation}
where $\var >0$ is the unbinding rate. 
\begin{prop}\label{prop-appli}
System \eqref{deux}-\eqref{cinq}, endowed by a continuous initial data $(u_0,v_0)$ satisfying \eqref{un} and \eqref{six}, and with $F(v)$ and $G(v)$ respectively defined by \eqref{Fdef} and \eqref{Gdef} has a unique homogeneous stable solution at equilibrium given by the couple 
$(U,V)$, with $U\geq 0$ defined by \eqref{UV} and $V\in ]0,\rho[$.
\end{prop}
{\bf Proof :}\\
It is clear that both $F$ and $G$ satisfy hypothesis a) to d). Hence the existence of the solution and relation \eqref{un} are ensured by Theorem \ref{main-theorem} and Theorem \ref{positivity}.  Moreover, the existence of the stationary homogeneous solution is proved in Proposition \ref{stationary-prop}. \\
Concerning the stability of the stationary solution, from Proposition \ref{stability-prop} we just have to prove \eqref{relation-h(v)}. Since $G'(v) = 0$, dividing \eqref{relation-h(v)} by $G(V)>0$, we obtain the stability of $(U,V)$ if the condition $F'(V) >0$ holds. And this is the case, since $F'(v)=b>0$ for all $v\in ]0,\rho[$. Concluding the proof. \ep

\begin{figure}[H]
\centering{\includegraphics[width=0.45\textwidth]{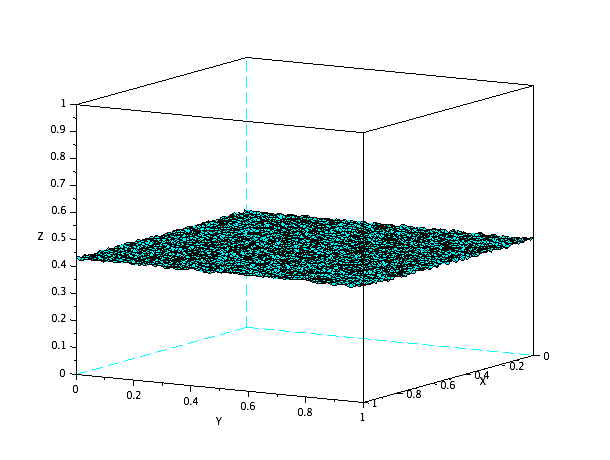}\quad
\includegraphics[width=0.45\textwidth]{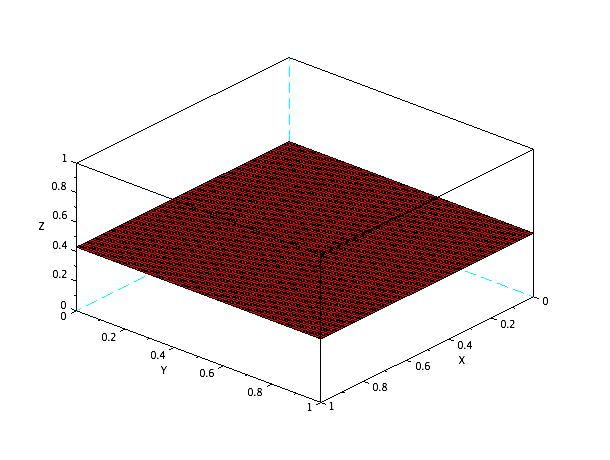}}
\caption{The initial data (left) and the equilibrium solution (right) for the function $v$. The parameters used for the computation are :
$\rho=1$, $a=0.25$, $b=0.5$ and  $\varepsilon=0.35$, $\sigma=1$ and $T=200$.}
\label{fig-numeric2}
\end{figure}

To conclude, we perform some numerical simulations defining the aggregation function $F(v)$ as in \eqref{Fdef} and the unbinding function as in \eqref{Gdef}. In Figure \ref{fig-numeric2} we show the convergence of the distribution function $v$ to the homogeneous function $v(x)=V$ when the initial data is defined as a small perturbation of the stationary solution $V$: 
$$v_0=V \pm \tilde{v},$$
with $|\tilde v|\approx 10^{-3}$. With our choice of parameters the stationary value is $V=0.43086$.

\begin{figure}[H]
\centering{\includegraphics[width=0.45\textwidth]{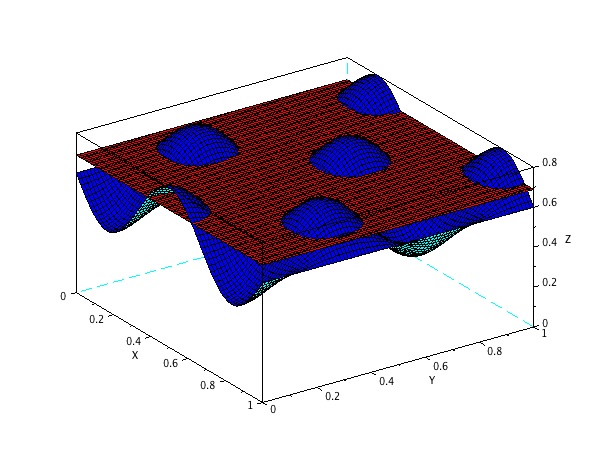}\quad
\includegraphics[width=0.45\textwidth]{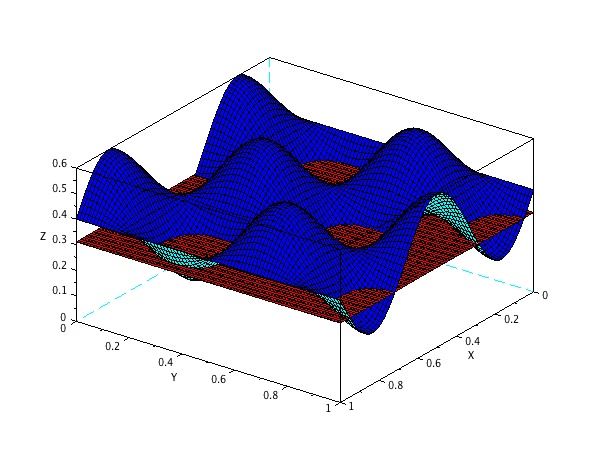}}
\caption{The initial data (blue) and the equilibrium solution (red) for the functions $u$ (left) and $v$ (right). The parameters used for the computation are :
$\rho=0.7$, $a=0.25$, $b=0.5$ and  $\varepsilon=0.35$, $\sigma=1$ and $T=10$.}
\label{fig-numeric}
\end{figure}

The evolution of the solution $(u,v)$ towards the stable equilibrium state for a 
$C^1(\Omega)$ initial data $(u_0,v_0)$ is shown in Figure \ref{fig-numeric}. We choose here a sinusoidal initial data defined by :
$$u_0(x)=1-v_0(x), \quad 
v_0(x) = 0.4+(0.2\ \sin(3 \pi x)\ \cos(3\pi y)).$$
Both particle distributions evolves to a homogeneous stationary state, which density value can be computed numerically solving $h(v)+v-1=0$ and is given, for our choice of parameters, by  $V=0.3107435$ and $U=0.6892565$.

\begin{figure}[H]
\centering{\includegraphics[width=0.3\textwidth]{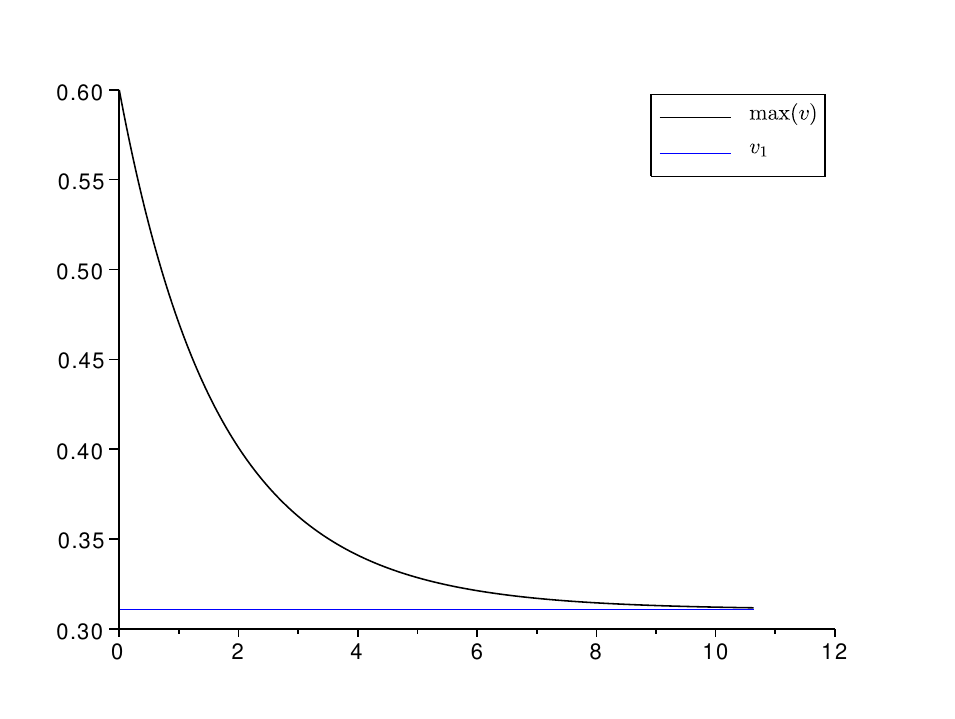}\
\includegraphics[width=0.3\textwidth]{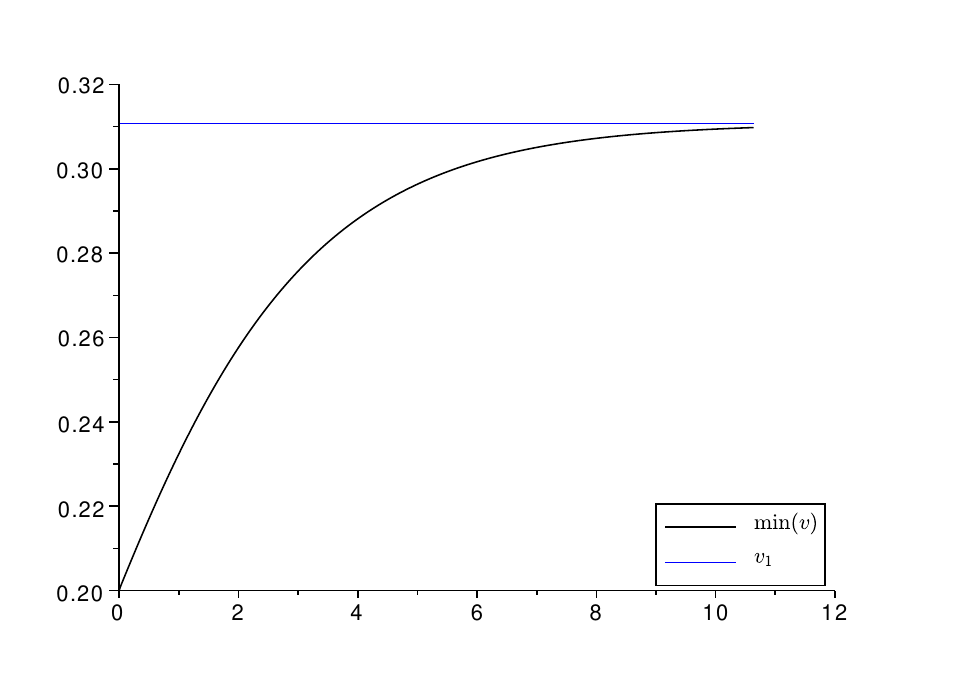}\
\includegraphics[width=0.3\textwidth]{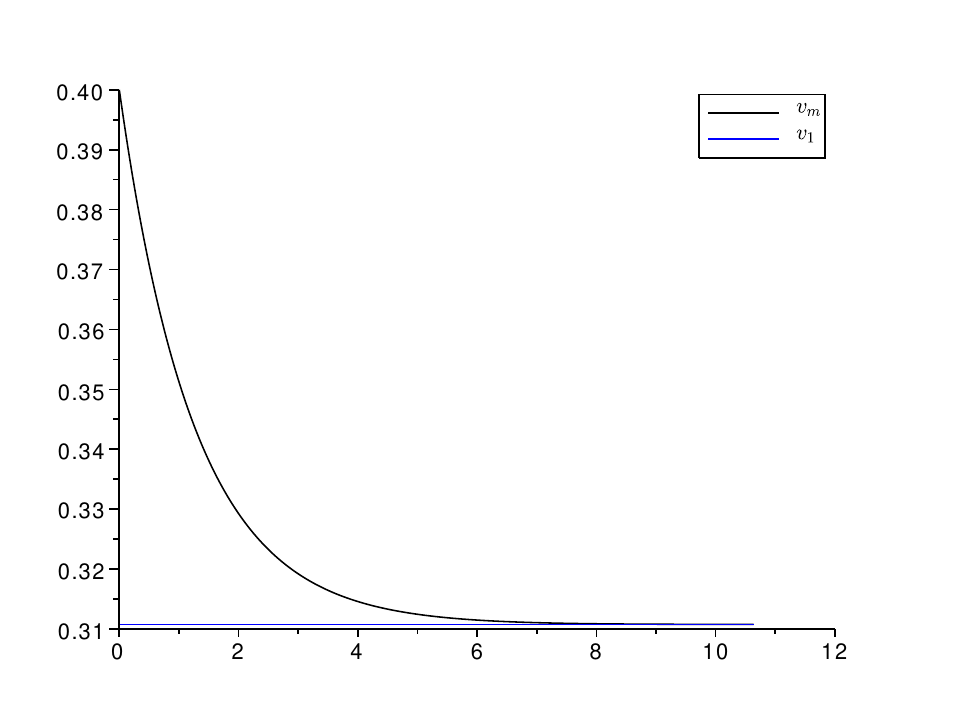}}
\caption{Convergence in time (black line) of the maximum (left), the minimum (middle) of $v(t,\cdot)$, and of mean value $v_m(t)$ (right) towards the value $v_1=0.3107435$ (blue line) for the parameters $\rho=0.7$, $a=0.25$, $b=0.5$,  $\varepsilon=0.35$ and $\sigma=1$.}\label{convergence}
\end{figure}

Finally, we may compute the convergence rate towards the stable solution. Figure \ref{convergence} shows the time convergence of the $\max(v)$, $\min(v)$ and of the mean value $v_m(t)$ towards $V=v_1$, . The computation has been stopped when the  maximum between the absolute errors $|v_0-\max(v)|$, $|v_0-\min(v)|$ and $|v_0-v_m|$ is smaller than $10^{-3}$. A linear regression study of the three black curves  shows that the convergences are exponential with a standard deviation of the order of $10^{-2}$.

This result shows that this simple local model of aggregation can't produce clusters as they are shown by experiments. Thus this choice is not suitable for the description of the biological problem we consider. Non-linear and local aggregation and unbinding rates leading to patterns formations are considered in \cite{SMSM} and shows very good agreements with experimental data.
  
\section*{Acknowledgments}
This work was supported by the PEPS-MBI project "Mod\'elisation d'Adh\'esion des Cadh\'erines" founded by the French research institutes CNRS, INRIA and INSERM, and by the Kibord project (ANR-13-BS01-0004) funded by the French Ministry of Research.

\end{document}